\documentclass{birkjour}

\usepackage{amsmath}
\usepackage{amssymb}
\usepackage{amsthm}
\usepackage{color}
\usepackage{graphicx}
\usepackage{psfrag}

\newcommand{\geqs}{\geqslant}
\newcommand{\leqs}{\leqslant}
\newcommand{\Nb}{{\mathbb N}}

\newenvironment{Proof}{\removelastskip \vskip12pt plus 1pt \noindent
{\em Proof.\/}\rm }{\hfill$\square$ \vskip12pt plus 1pt}

\newtheorem{theorem}{Theorem}[section]
\newtheorem{lemma}[theorem]{Lemma}

\theoremstyle{definition}
\theoremstyle{remark}

\numberwithin{equation}{section}

\begin{document}

%\title[The RBK cluster system: scaling behaviour]{The Redner--Ben-Avraham--Kahng coagulation system 
%with constant coefficients: the finite dimensional case}
\title[The finite dimensional RBK coagulation system]{The Redner--Ben-Avraham--Kahng coagulation system 
with constant coefficients: the finite dimensional case}%    Information for first author
\author{F.P. da Costa}
\address{Departamento de Ci\^encias e Tecnologia, Universidade Aberta, Lisboa, Portugal, and
Centro de An\'alise Matem\'atica, Geometria e Sistemas Din\^amicos, Instituto Superior T\'ecnico,
Universidade de Lisboa, Lisboa, Portugal}
%\curraddr{Department of Mathematics and Statistics,
%Case Western Reserve University, Cleveland, Ohio 43403}
\email{fcosta@uab.pt}
    
\thanks{We thank an anonymous referee for the careful reading of the manuscript and for 
calling our attention to a sloppy argument in the original proof of Lemma 3.3. 
This work was partially supported by FCT under Strategic Project - LA 9 - 2013-2014}


\author{J.T. Pinto}
\address{Departamento de Matem\'atica and
Centro de An\'alise Matem\'atica Geometria e Sistemas Din\^amicos,
Instituto Superior T\'ecnico, Universidade de Lisboa, Lisboa, Portugal}
\email{jpinto@math.tecnico.ulisboa.pt}

\author{R. Sasportes}
\address{Departamento de Ci\^encias e Tecnologia, Universidade Aberta, Lisboa, Portugal, and
Centro de An\'alise Matem\'atica, Geometria e Sistemas Din\^amicos, Instituto Superior T\'ecnico,
Universidade de Lisboa, Lisboa, Portugal}
\email{rafael@uab.pt}

%    General info

\date{January 15, 2014, and, in revised form, August 13, 2014.}


%=============================================================================
%                                                                                 Abstract
%=============================================================================

\begin{abstract}
We study the behaviour as $t\to\infty$ of solutions $(c_j(t))$ to the 
Redner--Ben-Avraham--Kahng coagulation system with positive and compactly supported initial data,
rigorously proving and slightly extending results originally established in \cite{RBK} by means
of formal arguments.
\end{abstract}

\maketitle

%==============================================================================
%                                                                                Introduction
%==============================================================================

\section{Introduction}

In a recent paper \cite{CPS} we started the study of a coagulation model first considered in \cite{war,RBK} which we have called the 
Redner--Ben-Avraham--Kahng cluster system (RBK for short). This is the infinite-dimensional ODE system 
\begin{equation}
 \frac{dc_j}{dt} = \sum_{k=1}^{\infty}a_{j+k,k}c_{j+k}c_k - \sum_{k=1}^{\infty}a_{j,k}c_{j}c_k, \qquad j=1,2,\dots.
\label{rbkgeneral}
\end{equation}
with symmetric positive coagulation coefficients $a_{j,k}$. As with the discrete Smoluchowski's coagulation system \cite{bc} 
this is a mean-field model describing the evolution 
of a system given at each instant by a sequence $(c_{j})$, such that $c_{j}$ is the density of $j$-clusters 
for each integer $j$,  undergoing a binary reaction  described by a bilinear infinite-dimensional vector field. However, while in
the Smoluchowski's coagulation model  one $k$-cluster reacts with one $j$-cluster producing one $(j+k)$-cluster, in RBK
the interaction between such clusters produce one $|k-j|$-cluster.

\begin{figure}[h]
\begin{center}
\psfrag{j}{$j$-cluster}
\psfrag{k}{$k$-cluster}
\psfrag{j-k}{$|j-k|$-cluster}
\psfrag{+}{$+$}
\includegraphics[scale=.75]{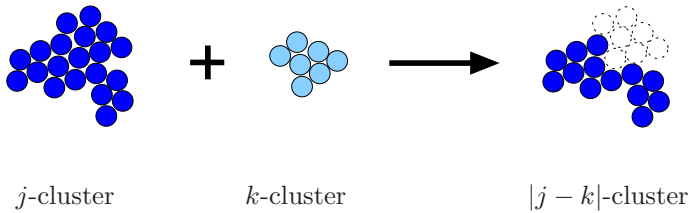}
\end{center}
\caption{Schematic reaction in the RBK coagulation model}\label{fig1}
\end{figure}

If we assume that there is no destruction of mass,
in the former model it makes sense to think of $j$ as the size, or mass, of each $j$-cluster. However in RBK the situation is 
different since with the same interpretation there would be a loss of mass in each reaction. Hence, it makes more sense to think of $j$ as
the size of the cluster `active part', being the difference between $(j+k)$ and $|j-k|$ the size of 
the resulting cluster that becomes inactive for the reaction process. 
A pictorial illustration of this is presented in Figure~\ref{fig1}.
%%%%%%%%%%%%%%%%%%%%%%%%%%%%%%%%%%%%%%%%%%%%%%%%%%%%%%%%%%%%%%%%%%%%%%%%%%%%%%%%%%%%
%
%
%
%
%%%%%%%%%%%%%%%%%%%%%%%%%%%%%%%%%%%%%%%%%%%%%%%%%%%%%%%%%%%%%%%%%%%%%%%%%%%%%%%%%%%%

For more on the physical interpretation of \eqref{rbkgeneral} see \cite{CPS, war, RBK}.

The nonexistence of a mass conservation property in RBK model makes for one of the major differences with respect to the 
Smoluchowski's model. Also, unlike in this one, in RBK a $j$ and a $k$-cluster react to produce 
a $j'$-cluster with $j'<\text{max}\{j,k\}$, 
implying that to an initial condition with an upper bound $N$ for the subscript values $j$ for which $c_j(0)>0$ there 
corresponds 
a solution with the same property for all instants $t\geqs 0$.    
This is an invariance property rigorously stated on Proposition 7.1 in \cite{CPS}. In this work we will consider such 
solutions for a finite prescribed upper bound $N\geqs 3$
and $j$-independent coagulation coefficients $a_{j,k}=1$, for all $j,k$. Then, if $c_{j}(0)=0$, for all $j\geqs N+1$, 
then $c_{j}(t)=0$ for $t\geqs 0$ and for the same values of $j$, while
$(c_{1}(t),c_{2}(t),\dots,c_{N}(t))$ satisfy the following  $N$-dimensional ODE
\begin{equation}
 \frac{dc_j}{dt} = \sum_{k=1}^{N-j}c_{j+k}c_k -c_{j} \sum_{k=1}^{N}c_k, \qquad j\in\Nb\cap [1,N], \label{rbk}
\end{equation}
where the first sum in the right-hand side is defined to be zero when $j=N$.

In this work we study system \eqref{rbk} for nonnegative initial conditions at $t=0$, from the point of view of the 
asymptotic behaviour of each component,
$c_{j}(t)$, $j=1,\dots,N$, as $t\to\infty$. This problem has already been addressed in \cite{RBK}, where the 
authors have used a formal approach.
 In Theorem \ref{theo},  we obtain the result for the general case $c_{j}(0)\geqs 0$, for $j=1,2,\dots,N$, 
proving rigorously that the result in  \cite{RBK} is correct for initial conditions such that $c_{N}(0)>0$ and 
the greater common divisor of the subscript values $j$ for which  $c_{j}(0)>0$ is $1$.  

%==============================================================================
%                                                                               The main result
%==============================================================================

\section{The main result}\label{main result}

Consider $N\geqs 3$. We are concerned with nonnegative solutions of \eqref{rbk}. By applying the results we have proved in \cite{CPS}
in the more general context refered above,  we can deduce that, 
for a solution $c=(c_j)$ to \eqref{rbk},
if $c_j(0)\geqs 0$, for $j=1,\dots,N$, then it is defined 
for all $t\in[0,\infty)$ and $c_j(t)\geqs 0$, for $j=1,\dots,N$, and all positive $t$.
Let $P=\{j\in\mathbb{N}\cap[1,N]\,\vert\,c_j(0)>0\}$ be   
the set of subscript values for which the components of the initial condition $c(0)$ are positive, and 
let $\operatorname{gcd}(P)$ be the greatest common divisor of the elements of $P$. In this paper we prove the following:

\begin{theorem}\label{theo}
Let $c=(c_j)$ be a solution of \eqref{rbk} satisfying $c_j(0)\geqs 0$ for all $j=1,\dots,N$. If $m:=\operatorname{gcd}(P)$ and
$p:=\sup P$, then,
for each $j=m,2m,\dots,p$, there exists 
$e_{j}:[0,\infty)\to \mathbb{R}$ such that $e_{j}(t)\to 0$ as $t\to \infty$, and
\[
c_{j}(t)=\frac{\widetilde{A}_{j}}{t(\log t)^{j/m-1}}(1+e_{j}(t))
\]
where
\[
\widetilde{A}_{j}:= \frac{(N-1)!}{(N-j/m)!}.
\]
For all other $j\in\mathbb{N}\cap [1,N]$, $c_j(t)=0$, for all $t\geqs 0$.
\end{theorem}

We begin the proof of this theorem by reducing it to the case $m=1$, $p=N$. Consider, for each $t\geqs 0$, 
$\mathcal{J}(t):=\{j\in\mathbb{N}\cap [1,N]\,\vert\, c_j(t)>0\}$, the set of subscript values for 
which the components of the solution are positive at instant $t$. Obviously, $P=\mathcal{J}(0)$. 
The case $\# P=1$ is an immediate consequence of Proposition 7.3 in \cite{CPS} and its proof. Consider now the case $\# P >1$.
Then, according to Proposition 7.2 in \cite{CPS},  $\mathcal{J}(t)=m\mathbb{N}\cap[1,p],$ for all $t>0$. Let $\tilde{N}:=p/m$ and, for
$j=1,2,\dots,\tilde{N}$, let us write $\tilde{c}_j:=c_{jm}$. Then it is straightforward to check that \eqref{rbk} 
is again satisfied with $N$ and $c_j$, for $j=1,2,\dots,N$, 
replaced by $\tilde{N}$ and $\tilde{c}_j$, for $j=1,2,\dots,\tilde{N}$, respectively. 
From the definition of $\mathcal{J}(t)$, we also have that, for $j=1,\dots,\tilde{N}$ and for all $t>0$, $\tilde{c}_j(t)>0$. 
For $j=1,\dots,N$, if $j\notin m\mathbb{N}\cap [1,p]$, then
$c_j(t)=0$, for all $t\geqs 0$.
Hence, after having established the validity of Theorem \ref{theo} with the restrictions $m=1$ and $p=N$, 
if we consider a solution $c(\cdot)$ with initial conditions for which $m>1$, $p<N$ or both,
we can apply that restricted version of the theorem to $\tilde{c}$ and then use the fact that, for $j=1,\dots,p$, $c_j(t)=\tilde{c}_{j/m}(t)$. 
For the other subscript values, $c_j(t)$ identically vanishes.

In conclusion, it is sufficient to prove the above theorem for $m=1$, $p=N$, in which case, as we have 
seen, $c_j(t)>0$, for $j=1,2,\dots,N$, and all $t>0$.  
This is done in next section.
%==============================================================================
%                                                                                Long time behaviour of strictly positive solutions
%==============================================================================

\section{Long time behaviour of strictly positive solutions}
\label{lemmas}
%Let us assume that $c_{j}(t)>0$, for all $t>0$ and all j=1,\dots,N. 
%According to Proposition 7.2 in \cite{CPS} this is true if $\text{gcd}(P)=1$, where
%$P:=\{j\in \mathbb{N}\,\vert\, c_{j}(0)>0\}$. This does not imply a lost of 
%generality since if $\text{gcd}(P)=m>1$, by the same Proposition, if we define
%$\tilde{c}_{j}=c_{mj}$, for $j=1,\dots, \tilde{N}:=N/m$,  then these new 
%variables satisfy \eqref{rbk} and $\tilde{c}_{j}(t)>0$, for all $t>0$. 
%For any $k$ not of the form $mj$, $c_{k}(t)=0,$ for all $t\geqs 0$.

Consider a solution $c(\cdot)=(c_j(\cdot))$ to
\eqref{rbk} such that $c_j(t)>0$ for all $j=1,\dots,N$ and all $t\geqs 0$. 
By the above and the fact that the ODE is autonomous we will see that this does not imply a loss of generality. 
Let
\[
\nu(t):=\sum_{j=1}^{N}c_{j}(t),
\]
so that \eqref{rbk} can be rewritten as
\begin{equation}\label{cjnu}
\dot{c}_{j}(t) +c_{j}(t)\nu(t) = \sum_{k=1}^{N-j}c_{j+k}(t) c_k(t),
\end{equation}
and, in particular,
\begin{equation}
\label{cN}
\dot{c}_{N}(t)+c_{N}(t)\nu(t)=0\,.
\end{equation}
We start by following the procedure already used in \cite{RBK} that consists in 
time rescaling \eqref{rbk} so that the resulting equations only retain the production terms. From \eqref{cN}
\[
c_{N}(t)/c_{N}(0)=\exp\left(-\int_{0}^t\nu(s)\,ds\right)\,.
\]
Since $e^{\int_{0}^t\nu}$ is an integrating factor of \eqref{cjnu}, we conclude that
\begin{equation}\label{cjcn}
\frac{d}{dt}\left(\frac{c_{j}(t)}{c_{N}(t)} \right)=\frac{1}{c_{N}(t)} \sum_{k=1}^{N-j}c_{j+k}(t) c_k(t)\,.
\end{equation}
Let $y(t):=\int_{0}^{t}c_{N}(s)\,ds$ and define functions $\phi_{j}(y)$, 
such that
\begin{equation}
 c_{j}(t)=\phi_{j}(y(t))c_{N}(t),    \label{phij}
\end{equation}
for each $j=1,\dots, N$, and $t\geqs 0$.
Then, for $j=1,\dots,N-1$, $\phi_{j}(y)$ is defined and is strictly positive for $y\in [0,\omega)$, where
$\omega:=\int_{0}^\infty c_{N}\in (0,+\infty]$.
Let us denote by $(\cdot)'$ the derivative with respect to $y$. Then, from \eqref{cjcn} we obtain
\begin{equation}
\label{phieq}
\begin{aligned}
\phi_{j}'(y)&= \sum_{k=1}^{N-j}\phi_{j+k}(y) \phi_k(y)\,,\quad j=1,\dots, N-1,\\
\phi_{N}(y)&=1\,,
\end{aligned}
\end{equation}
for $0\leqs y < \omega$. 
Conversely, if $(\phi_{j}(y))$ is a solution of \eqref{phieq} in 
its maximal positive interval $(0,\omega^{*})$ and if $c_{N}(\cdot)$, 
and therefore $y(\cdot)$, is given, then $c_{j}(t)=c_{N}(t)\phi_{j}(y(t))$, 
for $j=1,\dots, N$ solves \eqref{rbk} for $t\in[0,\infty)$, so that $\omega^{*}=\omega.$

In the next two lemmas we state some results about the asymptotic behaviour of $\phi(y)$. 
%
%===============================================================================
%                                             Lemma 1
%===============================================================================
%
\begin{lemma}\label{lem1}
Any solution of \eqref{phieq}, say $\phi(y)=(\phi_{1}(y),\dots,\phi_{N-1}(y),1)$, 
satisfying $\phi_{j}(0)>0$, for all $j=1,\dots,N$, is 
defined for $y\in[0,\omega)$ where $\omega>0$ is finite  and moreover,
\begin{enumerate}
\item[{\rm (i)}] $\phi_{j}(y)\to+\infty$ as $y\to\omega$, for all $j=1,2,\dots,N-1$;
\item[{\rm (ii)}] $\phi_{j}(y)/\phi_{j+1}(y)\to+\infty$ as $y\to\omega$, for all $j=1,2,\dots,N-1$.
\end{enumerate}
\end{lemma}
\begin{Proof}
Let $(\phi_{j}(y))$ be a solution of \eqref{phieq} in its positive maximal interval of 
existence $[0,\omega)$ satisfying the hypothesis of the lemma.
Then, for all $j=1,\dots, N$, $\phi_{j}(y)>0$, for all $y\in[0,\omega)$. Since, 
\begin{equation}\label{phieqapr}
\phi '_{j}(y)\geqs\phi_{j+1}(y)\phi_{1}(y)\,,
\end{equation}
for $j=1,\dots, N-1$ (with equality for $j=N-1$), and $\phi_{N}(y)=1$, by defining 
$\tau (y):=\int_{0}^{y}\phi_{1}(s)\,ds$, and $\psi_{j}(\tau),$ such that $\phi_{j}(y)=\psi_{j}(\tau(y))$, we obtain,
\begin{equation}\label{psieq}
\frac{d}{d\tau}\psi_{j}(\tau)\geqs \psi_{j+1}(\tau),
\end{equation}
for $j=1,\dots,N-1$ (with equality for $j=N-1$), $\psi_{N}(\tau)=1$, for $0\leqs \tau< \int_{0}^{\omega}\phi_{1}$. 
The $N-1$ equation gives,
\[
\psi_{N-1}(\tau)= \tau+c_{0}.
\]
Then by successively integrating \eqref{psieq} for $j=N-2,\,N-3,\dots,1$, and taking in account 
that $\psi_{j}(0)\geqs 0$ for $j=1,\dots,N$, we obtain
\[
\psi_{N-k}(\tau)\geqs \frac{\tau^k}{k!},\quad k=1,\dots, N-1\,.
\]
In particular,
\[
\psi_{1}(\tau)\geqs \frac{\tau^{N-1}}{(N-1)!}\,,
\]
which is equivalent to
\[
\tau'(y)\geqs \frac{ \tau(y)^{N-1}}{(N-1)!}\,.
\]
Since, by hypothesis, $N-1>1$, the last inequality means that $\tau(\cdot)$ blows up
at a finite value of $y$, which implies that $\omega<+\infty$. 
By fundamental results in ODE theory, this in turn implies that, for our solution, we have
$\|\phi(y)\|\to \infty$, as $y \to \omega$. This, together with the monotonicity 
property of each $\phi_{j}(y)$, implies that there is a $j^{*}\in\{1,\dots,N-1\}$ 
such that $\phi_{j^{*}}(y)\to +\infty$ as $y\to\omega$.
We now prove the nontrivial fact that this is true for all $j=1,\dots,N-1$. 
In order to derive such a conclusion we first prove that, for $j=1,\dots, N-1$, 
$\phi_{j}(y)/\phi_{j+1}(y)$ is bounded away from zero for $y$ sufficiently 
close to $\omega$. Specifically, we prove that for $n=N-1, N-2,\dots,2,1$, 
there are $\eta>0$, $Y\in[0,\omega)$ such that
\begin{equation}\label{phieta}
\frac{\phi_{j}(y)}{\phi_{j+1}(y)}>\eta,
\end{equation}
for $j=n,n+1,\dots, N-1$, and for all $y\in [Y,\omega)$.

Consider $n=N-1$. Then $\phi '_{N-1}(y)=\phi_{1}(y)$, so that $\phi_{N-1}(y)/
\phi_{N}(y)=\phi_{N-1}(0)+\int_{0}^{y}\phi_{1}$ and, by the positivity of $\phi_{1}$ 
the result is obvious with $\eta=\phi_{N-1}(Y)$ for any $Y\in (0,\omega)$.

Suppose now that we have proved our claim for  $n+1$, with $n\in\{1,\dots,N-1\}$, 
that is, there are $\eta>0$, $Y\in[0,\omega)$ such that  \eqref{phieta} is true, 
for $j=n+1,n+2,\dots, N-2$ and for  $y\in [Y,\omega)$. We prove the same holds for $n$. 
Since, for $y\in[Y,\omega)$
\[
\frac{\phi '_{n}(y)}{\phi '_{n+1}(y)}=\frac{ \displaystyle\sum_{k=1}^{N-n}
\phi_{k+n}(y) \phi_k(y)}{\displaystyle \sum_{k=1}^{N-n-1}\phi_{k+n+1}(y) \phi_k(y)}\geqs
\frac{ \displaystyle\sum_{k=1}^{N-n-1}\phi_{k+n+1}(y) \phi_k(y)\cdot 
\frac{\phi_{k+n}(y)}{\phi_{k+n+1}(y)}}{\displaystyle \sum_{k=1}^{N-n-1}\phi_{k+n+1}(y) \phi_k(y)}\geqs \eta,
\]
and therefore
\[
\phi '_{n}(y)\geq\eta\phi '_{n+1}(y),
\]
by integration we obtain
\[
\phi_{n}(y)-\phi_{n}(Y)\geqs \eta (\phi_{n+1}(y)-\phi_{n+1}(Y))
\]
or
\[
\frac{\phi_{n}(y)}{\phi_{n+1}(y)}\geqs \frac{\phi_{n}(Y)}{\phi_{n+1}(y)}+
\eta\left(1-\frac{\phi_{n+1}(Y)}{\phi_{n+1}(y)}\right)\,.
\]
Let $\tilde{Y}\in (Y,\omega)$. Then, for $y\in[\tilde{Y},\omega)$,
\[
\phi_{n+1}(y)\geqs \phi_{n+1}(\tilde{Y})>\phi_{n+1}(Y),
\]
and defining \[\tilde{\eta}:=\eta\left(1-\frac{\phi_{n+1}(Y)}{\phi_{n+1}(\tilde{Y})}\right)\]
we conclude that, for $y\in [\tilde{Y},\omega)$,
\[
\frac{\phi_{n}(y)}{\phi_{n+1}(y)}\geqs\tilde{\eta}.
\]
By redefining $Y,\eta$ as $\tilde{Y},\tilde{\eta}$ we have proved \eqref{phieta} for $n$. This completes our induction argument.

Now let $K:=\{j=1,\dots, N-1\,\vert\; \phi_{j}(y)\to\infty\text{ as }y\to\omega\}$. 
We already know that $K\not=\emptyset$, so that we can define
$J:=\max K$. Then, from \eqref{phieta} we get
\[
\phi_{j}(y)\to\infty\text{ as }y\to\omega,\quad \text{for all}\; j=1,\dots,J\,.
\]
It is then sufficient to prove that, in fact, $J=N-1$. This is based on the integral version of \eqref{phieq}, namely
\begin{multline}\label{phiint}
\phi_{j}(y)-\phi_{j}(Y)=\int_{Y}^{y}\phi_{j+1}\phi_{1}+\int_{Y}^{y}\phi_{j+2}
\phi_{2}+\dots\\+\int_{Y}^{y}\phi_{N-j-1}\phi_{N-1}+\int_{Y}^{y}\phi_{N-j},
\end{multline}
for $j=1,\dots,N-1$. Now, in order to derive a contradiction, suppose that $J<N-1$. 
Then, for $j=J+1,\dots,N-1$, $\phi_{j}(y)$ is bounded for $y\in[Y,\omega)$. 
But then, since \eqref{phiint} implies that
\begin{equation}\label{ineqast}
\phi_{j}(y)-\phi_{j}(Y)> \int_{Y}^{y}\phi_{N-j},
\end{equation}
we conclude that $\int_{Y}^{y}\phi_{j}$ must be bounded for $j=1,2,\dots,N-J-1$ 
and $y\in[Y,\omega)$. Therefore, by the monotonicity of all the $\phi_j(\cdot)$,
we get, for all $y\in [Y,\omega)$,
\[\begin{aligned}
\phi_{J}(y)-\phi_{J}(Y)&\leqs \phi_{J+1}(y)\int_{Y}^{y}\phi_{1}+\phi_{J+2}(y)\int_{Y}^{y}\phi_{2}+\dots\\
&\qquad\qquad\qquad\qquad\ldots+\phi_{N-1}(y)\int_{Y}^{y}\phi_{N-J-1}+\int_{Y}^{y}\phi_{N-J}\\
&\leqs M+\int_{Y}^{y}\phi_{N-J}, &
\end{aligned}\]
for some positive constant $M$. Since $\phi_{J}(y)\to \infty$, as $y\to\omega$, this 
bound forces $\int_{Y}^y\phi_{N-J}\to \infty$ as $y\to\omega$. Now,
again by \eqref{phieta}, we have, for $y\in [Y,\omega)$,
\[
\phi_{1}(y)\geqs \eta \phi_{2}(y)\geqs \eta^2 \phi_{3}(y)\geqs\dots\geqs\eta^{N-J-1}\phi_{N-J}(y)\,,
\]
implying that, for all $j=1, 2, \ldots, N-J-1,$
\[
\int_{Y}^y\phi_{j}\geqs \eta^{N-J-j}\int_{Y}^y\phi_{N-J},
\]
contradicting the boundedness conclusion following inequality \eqref{ineqast}. This proves that $J=N-1$.

It remains to be proved assertion {\rm (ii)}. For $j=N-1$ it is trivial, since
\[
\frac{\phi_{N-1}(y)}{\phi_{N}(y)}=\phi_{N-1}(y)\to+\infty\quad\text{as}\quad y\to\omega,
\]
as we have seen before. Suppose we have proved {\rm (ii)} for $j=N-1,N-2,\dots,n+1$ for some $n\in\{1,2,\dots,N-2\}$. 
We prove that the same holds for $j=n$.
 We consider again, for $y$ close to $\omega$, the quotient
\[\begin{aligned}
\frac{\phi '_{n}(y)}{\phi '_{n+1}(y)}&=\frac{ \displaystyle\sum_{k=1}^{N-n}
\phi_{k+n}(y) \phi_k(y)}{\displaystyle \sum_{k=1}^{N-n-1}\phi_{k+n+1}(y) \phi_k(y)}=
\frac{ \displaystyle\sum_{k=1}^{N-n}\frac{\phi_{k+n}(y)}{\phi_{2+n}(y)}\cdot
\frac{\phi_k(y)}{\phi_{1}(y)}}{1+\displaystyle \sum_{k=2}^{N-n-1}
\frac{\phi_{k+n+1}(y)}{\phi_{2+n}(y)}\cdot\frac{\phi_k(y)}{\phi_{1}(y)}}\\
&>\frac{\phi_{1+n}(y)}{\phi_{2+n}(y)}\left(1+\sum_{k=2}^{N-n-1}\eta^{-k+1}
\frac{\phi_{k+n+1}(y)}{\phi_{2+n}(y)}\right)^{-1}\to+\infty,
\end{aligned}\]
as $y\to\omega$. Then, we know by Cauchy's rule that
\[
\lim_{y\to \omega}\frac{\phi_{n}(y)}{\phi_{n+1}(y)}=\lim_{y\to\omega}
\frac{\phi '_{n}(y)}{\phi '_{n+1}(y)}=+\infty,
\]
and our induction argument is complete.
\end{Proof}
%
%===============================================================================
%
\begin{lemma}\label{lem2}
In the conditions of the previous lemma, for each $j=1,\dots,N-1$, there is 
$\rho_{j}:[0,\omega)\to \mathbb{R}$ such that $\rho_{j}(y)\to 0$ as $y\to \omega$, and
\[
\phi_{j}(y)=\frac{A_{j}}{(\omega-y)^{\alpha_{j}}}(1+\rho_{j}(y))\,,
\]
where
\[
\alpha_{j}:=\frac{N-j}{N-2},\qquad
A_{j}:=\frac{1}{(N-j)!}\left(\frac{(N-1)!}{N-2}\right)^{\alpha_{j}}\,.
\]
\end{lemma}
\begin{Proof}
By {\rm (ii)} of the previous lemma, we know that, for $j=1,\dots,N-1$,
\[
\frac{\displaystyle\sum_{k=1}^{N-j}\phi_{j+k}(y) \phi_k(y)}{\phi_{j+1}(y)\phi_{1}(y)}=
1+\sum_{k=2}^{N-j}\frac{\phi_{j+k}(y)}{\phi_{j+1}(y)}\cdot\frac{\phi_{k}(y)}{\phi_{1}(y)}\to 1\quad\text{as}\quad y\to\omega\,.
\]
Hence, we can write, for $j=1,\dots,N-1$, and $y\in(0,\omega)$
\begin{equation}\label{phieqr}
\phi '_{j}(y)=\phi_{1+j}(y)\phi_{1}(y)(1+r_{j}(y))
\end{equation}
such that $r_{j}(y)\to 0$, as $y\to\omega$. We now perform the same change of variables as in the 
beginning of the proof of the previous lemma, this time giving, for $\tau\geqs 0$,
\begin{equation}\label{psieqr}
\frac{d}{d\tau}\psi_{j}(\tau)= \psi_{j+1}(\tau)(1+\hat{r}_{j}(\tau)),
\end{equation}
such that $\hat{r}_{j}(\tau)\to 0$, as $\tau\to\infty$.  We now prove that, for $j=1,\dots,N-1$, 
\begin{equation}\label{psirho}
\psi_{j}(\tau)=\frac{\tau^{N-j}}{(N-j)!}(1+\hat\rho_{j}(\tau))
\end{equation}
where $\hat\rho_{j}(\tau)\to 0$ as $\tau\to\infty$. For $j=N-1$, taking into 
account that $\hat{r}_{N-1}(\tau)\equiv 0$, 
the result easily follows:
\[
\psi_{N-1}(\tau)=\tau+c_{0}=\tau(1+c_{0}\tau^{-1})\,.
\]
Now suppose we have verified \eqref{psirho} for $j=n+1$, for some $n=1,\dots,N-2$. 
We prove the same holds for $j=n$.
Defining $\delta (\tau)$ by
\[
\delta(\tau)=(1+\hat{\rho}_{n+1}(\tau))(1+\hat{r}_{n}(\tau))-1\,,
\]
we have $\delta(\tau)\to 0$ as $\tau\to \infty$, and by \eqref{psieqr} and \eqref{psirho},
\[
\frac{d}{d\tau}\psi_{n}(\tau)=\frac{\tau^{N-n-1}}{(N-n-1)!}(1+\delta(\tau))\,,
\]
and therefore, upon integration,
\[
\psi_{n}(\tau)-\psi_{n}(0)=\frac{\tau^{N-n}}{(N-n)!}+\frac{1}{(N-n-1)!}\int_{0}^{\tau}s^{N-n-1}\delta(s)\,ds,
\]
which can be written as
\[
\psi_{n}(\tau)=\frac{\tau^{N-n}}{(N-n)!}(1+\hat{\rho}_{n}(\tau))\
\]
where
\[
\hat{\rho}_{n}(\tau):=\frac{(N-n)!\psi_{n}(0)}{\tau^{N-n}}+\frac{N-n}
{\tau^{N-n}}\int_{0}^{\tau}s^{N-n-1}\delta(s)\,ds\,.
\]
If the integral in the right hand side stays bounded for $\tau\geqs 0$, then the last term converges to $0$ as 
$\tau\to\infty$. If it is unbounded, since its integrand is positive then the integral tends to $+\infty$, as 
$\tau\to\infty$. In this case we can apply Cauchy's rule since
\[
\frac{\left(\int_{0}^{\tau}s^{N-n-1}\delta(s)\,ds\right)'}{(\tau^{N-n})'}=\frac{\delta(\tau)}{N-n}\to 0,
\quad\text{as}\quad \tau\to\infty\,,
\] 
thus proving that also in this case, the last term converges to $0$ as $\tau\to\infty$. Either way we have 
$\hat{\rho}_{n}(\tau)\to 0$ as $\tau\to\infty$, thus proving assertion \eqref{psirho} for $j=n$. 
Our induction argument is complete.

In particular, 
\[
\psi_{1}(\tau)=\frac{\tau^{N-1}}{(N-1)!}(1+\hat\rho_{1}(\tau))
\]
which is equivalent to
\[
\tau '(y)=\frac{\tau(y)^{N-1}}{(N-1)!}(1+\hat\rho_{1}(\tau(y)))
\]
for $y\in (0,\omega)$.

Let $0<y<y_{1}< \omega$. Then, the integration of the previous equality in $[y,y_{1}]$ yields
\[
\tau(y)^{2-N}-\tau(y_{1})^{2-N}=\frac{N-2}{(N-1)!}\left( y_{1}-y+\int_{y}^{y_{1}}\hat{\rho}_{1}(\tau(s))\,ds\right)\,.
\]
Define  $\hat{R}(y,y_{1}):=\frac{1}{y_{1}-y}\int_{y}^{y_{1}}\hat{\rho}_{1}(\tau(s))\,ds$. Then,
\begin{equation}\label{intphi1ineq}
\tau(y)=\left[\tau(y_{1})^{2-N}+\frac{N-2}{(N-1)!}(y_{1}-y)(1+\hat{R}(y,y_{1}))\right]^{-\frac{1}{N-2}}\,.
\end{equation}
Now, observe that $\tau(y_{1})^{2-N}\to 0,$ as $y_{1}\to\omega$. Also, by fixing $y\in (0,\omega)$, 
for $y_{1}\in [y+\eta,\omega)$ with $\eta>0$ small, $y_{1}\mapsto \hat{R}(y,y_{1})$ is bounded.
Therefore we can define $R_{0}(y):=\lim_{y_{1}\to\omega}\hat{R}(y,y_{1})$.  Then by making $y_{1}\to\omega$ 
in \eqref{intphi1ineq} we obtain
\begin{equation}\label{tauy}
\tau(y)=\left[\frac{N-2}{(N-1)!}(\omega-y)(1+R_{0}(y))\right]^{-\frac{1}{N-2}}\,.
\end{equation}
with
\[
R_{0}(y)=\frac{1}{\omega-y}\int_{y}^{\omega}\hat{\rho}_{1}(\tau(s))\,ds \to 0\quad\text{as}\quad y\to \omega,
\]
by Cauchy rule and the fact that $\hat{\rho}_{1}(\tau(y))\to 0$ as $y\to\omega$.

For $j=1,\dots, N-1$, define
\[
\rho_{j}(y):=\left(1+R_{0}(y)\right)^{-\frac{N-j}{N-2}}\left(1+\hat{\rho}_{j}(\tau(y))\right)-1\,.
\]
so that $\rho_{j}(y)\to 0$, as $y\to \omega$. By \eqref{psirho} and \eqref{tauy}, for $j=1,\dots, N-1$ and $y\in(0,\omega)$,
\[
\phi_{j}(y)=\psi_{j}(\tau(y))
=\frac{1}{(N-j)!}\left(\frac{(N-1)!}{N-2}\right)^{\frac{N-j}{N-2}}(\omega-y)^{-\frac{N-j}{N-2}}(1+\rho_{j}(y))\,
\]
and the proof is complete.
\end{Proof}

The following lemma is a weaker version of Theorem \ref{theo} which will be used to complete the proof of the full result:

%%%%%%%%%%%%%%%%%%%%%%%%%%%%%%%%%%%%%%%%%%%%%%%%%%%%%%%%%%%%%%%%%%%%%%%%%%%%%%%%%%%%%%%%%%%%%%%
\begin{lemma}\label{lem3}
If $c_{j}(0)>0$, for $j=1,\dots,N$, then, for each such $j$, there exists 
$e_{j}:[0,\infty)\to \mathbb{R}$ such that $e_{j}(t)\to 0$ as $t\to \infty$, and
\[
c_{j}(t)=\frac{\widetilde{A}_{j}}{t(\log t)^{j-1}}(1+e_{j}(t))
\]
where
\[
\widetilde{A}_{j}:= \frac{(N-1)!}{(N-j)!}.
\]
\end{lemma}

\begin{Proof}
It was proved in \cite{CPS} that 
\[
 \nu_{\text{odd}}(t):=\sum_{\substack{ j=1\\j \text{\,odd}}}^{N}c_j(t)
\]
satisfies the differential equation $\dot{\nu}_{\text{odd}}=-\nu_{\text{odd}}^2$, and thus
\[
\nu_{\text{odd}}(t)=\frac{1}{(\nu_{\text{odd}}(0))^{-1}+t}.
\]
It follows that
\[
\nu_{\text{odd}}(t)=\frac{1}{t}(1+o(1))\qquad\text{as } t\to\infty.
\]
Defining $\nu_{\text{even}}(t)=\sum_{j=2,j \text{ even}}^{N}c_{j}(t)$ and using  Lemma~\ref{lem1}(ii) we have
\[
\frac{\nu_{\text{even}}(t)}{\nu_{\text{odd}}(t)}=\frac{\frac{c_{2}}{c_{1}}+\frac{c_{4}}{c_{1}}+\cdots + 
\frac{c_{2\lfloor N/2\rfloor}}{c_{1}}}
{1+\frac{c_{3}}{c_{1}}+\cdots+ \frac{c_{\mbox{}_{2\lfloor (N-1)/2\rfloor+1}}}{c_{1}}}=o(1), \qquad\text{as } t\to\infty.
\]
It follows that,
as $t\to\infty,$ 
\begin{equation}
\nu (t)=\nu_{\text{odd}}(t)\left(1+\frac{\nu_{\text{even}}(t)}{\nu_{\text{odd}}(t)}\right)=\nu_{\text{odd}}(t)(1+o(1))=\frac{1}{t}(1+o(1)).
\label{nu1}
\end{equation}
On the other hand, again by Lemma~\ref{lem1}(ii) and \eqref{phij}, we conclude that,  as $t\to\infty,$
\begin{equation}
\nu(t)=\sum_{j=1}^{N}c_{j}(t)=c_{1}(t)\left(1+\sum_{j=2}^{N}\frac{c_{j}(t)}{c_{1}(t)}\right)=c_{1}(t)(1+o(1)).
\label{nu2}
\end{equation}
From \eqref{nu1} and \eqref{nu2} we conclude that 
\[
tc_{1}(t)\to 1, \qquad\text{as } t\to\infty.
\]
By \eqref{phij} with $j=1,$ we can write $c_{1}(t)=\phi_{1}(y(t))c_{N}(t),$ and thus 
\begin{equation}
t\phi_{1}(y(t))c_{N}(t)\to 1, \qquad\text{as } t\to\infty.
\label{tphi1cN}
\end{equation}
When $j=1$, Lemma~\ref{lem2} reduces to
\begin{equation}
\phi_{1}(y)=\frac{A_{1}}{(\omega -y)^{\frac{N-1}{N-2}}}(1+o(1)), \qquad\text{as } y\to\omega.
\label{phi1}
\end{equation}
From  \eqref{tauy} we have
\(
\omega -y=\frac{(N-1)!}{N-2}\tau(y)^{2-N}(1+o(1)), \text{ as } y\to\omega,
\)
where $\tau(y)$ was defined by $\tau(y)=\int_{0}^{y}\phi_{1}(\tilde{y})d\tilde{y}$ in the beginning of the 
proof of Lemma~\ref{lem1}, and hence
\[
\tau(y(t))=\int_{0}^{y(t)}\phi_{1}(\tilde{y})d\tilde{y}=\int_{0}^{t}\phi_{1}(y(s))c_N(s)ds=\int_{0}^{t}c_{1}(s)ds.
\]
Since
\[
\frac{(\tau(y(t)))'}{(\log t)'}=\frac{c_{1}(t)}{1/t}=tc_{1}(t)\to 1, \qquad\text{as } t\to\infty,
\]
using Cauchy's rule we have
\(
\tau(y(t))=(\log t)(1+o(1)), \text{ as } t\to\infty,
\)
so that 
\begin{equation}
\omega -y(t)=\frac{(N-1)!}{N-2}(\log t)^{2-N}(1+o(1)), \qquad\text{as } t\to\infty,
\label{omy}
\end{equation}
and by \eqref{phi1}
\[
\phi_{1}(y(t))=A_{1}\left(\frac{N-2}{(N-1)!}\right)^{\frac{N-1}{N-2}}(\log t)^{N-1}(1+o(1)), \qquad\text{as } t\to\infty.
\]
Multiplying by $tc_{N}(t)$ and recalling \eqref{tphi1cN} we have
\[
A_{1}\left(\frac{N-2}{(N-1)!}\right)^{\frac{N-1}{N-2}}(\log t)^{N-1}tc_{N}(t)(1+o(1))\to 1, \qquad\text{as } t\to\infty,
\]
and since $A_{1}\left(\frac{N-2}{(N-1)!}\right)^{\frac{N-1}{N-2}}=\frac{1}{(N-1)!}$, we obtain
\[
\frac{t(\log t)^{N-1}}{(N-1)!}c_{N}(t)(1+o(1))\to 1, \qquad\text{as } t\to\infty,
\]
and it follows that,
as $t\to\infty,$ 
\begin{equation}
 c_N(t)  =  ((N-1)!)\frac{1}{t(\log t)^{N-1}}(1+o(1)).
\label{cn}
\end{equation}
Now we can use \eqref{phij}, Lemma~\ref{lem2}, and \eqref{cn} to obtain
\[
c_j(t) =
 \frac{A_{j}}{(\omega-y(t))^{\alpha_{j}}}
{\displaystyle((N-1)!)\,\frac{1}{t(\log t)^{N-1}}(1+o(1))}\qquad\text{as $t\to\infty$,}
\]
and from this, using \eqref{omy} and the definitions of $\alpha_j$ and $A_j$ in the statement of Lemma~\ref{lem2}, it follows that
\begin{equation}
 c_j(t) = \frac{(N-1)!}{(N-j)!}\,\frac{1}{t(\log t)^{j-1}}(1+o(1))\qquad\text{as $t\to\infty$,}
\end{equation}
as we wanted to prove.
\end{Proof}

Now, consider the case $c_{j}(0)\geqs 0$, for $j=1,\dots, N$, with $m=\operatorname{gcd}(P)=1$ and $p=\sup P=N$, thus implying that
$\mathcal{J}(t)=\mathbb{N}\cap [1,p]$ for all $t>0$. Since \eqref{rbk} is an autonomous ODE, then, given a small $\varepsilon>0$, for $t\geqs \varepsilon$, $c(t)=c_{\varepsilon}(t-\varepsilon)$,
where $c_{\varepsilon}(\cdot)$ is the solution of \eqref{rbk} satisfying the initial condition $c_{\varepsilon}(0)=c(\varepsilon)$. Therefore, the conditions of Lemma \ref{lem3} apply to $c_{\varepsilon}(\cdot)$. Then, it is easy to see that the asymptotic results that we conclude with respect to $c_{\varepsilon}(t)$ also apply to $c(t)$, allowing us to state the following:
\begin{lemma}
Let $c=(c_{j})$ be a solution satisfying $c_{j}(0)\geqs 0$, with $m=1$ and $p=N$. Then the conclusions of Lemma~$\ref{lem3}$ hold.
\end{lemma} 

This is, in fact, the particular case of Theorem \ref{theo} from which the full case follows as stated at the end of section \ref{main result}.  
%
%====================================================================================
%                                                                           Final remarks
%====================================================================================
%
\section{Final remarks}
A natural question to ask is: what is the asymptotic behaviour of the solutions of 
\eqref{rbk} in the infinite dimensional case ($N=\infty$)?
It is clear that Theorem \ref{theo} by itself is unsufficient to answer this question 
since the passage to the limit, $N\to \infty$, is not allowed without results on the 
uniformity of the various limits involved, which seems to be a hard task. Also it is 
far from clear how to rebuild the proofs of the lemmas in section \ref{lemmas} in this 
more general case since they heavily rely on the fact that there is a `last equation', 
the $N$-component equation, that  can be integrated by the reduction method we have used, 
being the asymptotic behaviour of the other components deduced in a `backwards' manner. 
Such procedure is obviously impossible in an infinite dimensional setting. In fact, that 
the situation can be very different for $N=\infty$ from the one displayed by Theorem 
\ref{theo} is shown by the existence of the self-similar solutions given by,
\[
c_{j}(t)=(\kappa+t)^{-1}(1-\alpha^2)\alpha^{j-1},\quad j=1,2,\dots,\quad t\geqs 0,
\]
with constants $\kappa>0$ and $\alpha\in(0,1)$ (see \cite{CPS}), in which case, 
$tc_{j}(t)\to (1-\alpha^2)\alpha^{j-1}$, as $t\to\infty$, for $j=1,2,\dots$. 
Further work will be devoted to fully understand this problem.

%
%====================================================================================
%                                                                           References
%====================================================================================
%
\bibliographystyle{amsplain}

\begin{thebibliography}{10}


\bibitem{bc} J.M. Ball, J. Carr, \textit{The discrete coagulation-fragmentation equations: existence, 
uniqueness, and density conservation}, J. Stat. Phys. {\bf 61}, 1/2 (1990) 203--234.


\bibitem{CPS} F.P. da Costa, J.T. Pinto, R. Sasportes, \textit{The Redner--Ben-Avraham--Kahng cluster system},
S\~ao Paulo J. Math. Sci. \textbf{6}, 2 (2012) 171--201.

\bibitem{war} I. Ispolatov, P.L. Krapivsky, S. Redner, \textit{War: The dynamics of vicious civilizations},
Phys. Rev. E \textbf{54} (1996) 1274--1289.

\bibitem{RBK} S. Redner, D. Ben-Avraham, B. Kahng, \textit{Kinetics of `cluster eating'}, J. Phys. A: Math. Gen.
\textbf{20} (1987), 1231--1238.


\end{thebibliography}

\end{document}